\documentclass[12pt]{article}
\usepackage{amsmath,amsfonts}
\begin{document}

\baselineskip 20pt

\title{ A proof of Fatou's interpolation theorem}

\author{Arthur A.~Danielyan}

\maketitle

\begin{abstract}

\noindent For Fatou's interpolation theorem of 1906 we suggest a new elementary proof.

\end{abstract}

\begin{section}{Introduction}

Denote by $\Delta$ and $T$ the open unit disk and the unit circle in the complex plane,
respectively. Recall that the disk algebra $A$ is the algebra of all functions
 on the closed unit disk $\overline \Delta$ that are
analytic on $\Delta$. Let $m$ be the Lebesgue measure on $T$.

The following theorem of Fatou \cite{Fa} is a cornerstone in many classical and contemporary
 investigations in complex analysis of one and several 
 variables as well as in functional analysis.

 \vspace{0.25 cm}
  
{\bf Theorem A.}\ {\it Let $F$ be a closed subset of $T$ such that $m(F)=0$. 
Then there exists a function $\omega_F$ in the disk algebra $A$ which vanishes precisely on $F$.}

\vspace{0.25 cm}

A slightly alternative formulation of Theorem A is the following: 

 \vspace{0.25 cm}  
  
{\bf Theorem A$'$.}\ {\it Let $F$ be a closed subset of $T$ such that $m(F)=0$. 
Then there exists a function $\lambda_F$ in the disk algebra $A$ such that $\lambda_F(z)=1$ on $F$ and $|\lambda_F(z)|<1$ on $T \setminus F$.}

\vspace{0.25 cm}

An early applications of Fatou's theorem is the original proof of the F. and M. Riesz theorem 
on the analytic measures, which asserts that if $\mu$ is a complex 
regular Borel measure on the unit circle orthogonal to all polynomials (in $z$),
then $\mu$ is absolutely continuous with respect to the Lebesgue measure. 
In fact the proof of the Riesz brothers is
the simplest among the other proofs of their theorem (cf. \cite{Koos}, p. 28).

Fatou's proof of Theorem A still remains the only proof of the theorem. It
has been included in many books (such as \cite{Gar}, \cite{Hoff}, \cite{Koos}, \cite{Priv} to list a few).
For the proof, Fatou \cite{Fa} first constructs  a summable function which is 
continuously differentiable in the complement of $F$ and continuous in the extended sense at the points of $F$.
This construction actually is the crucial step of the proof as the final part of the proof is standard.

In this paper we present a new proof for Theorem A. 
The proof is based on a summable function which we construct using a new method.
To complete the proof of Theorem A, we borrow the classical argument
on the existence of the boundary values of the conjugate harmonic functions, 
but, in general, the present proof
is both shorter and simpler than the original proof of Fatou.

\end{section}

\begin{section}{The proof}

The main idea of our proof is the following proposition, the proof of which is obvious. 

\vspace{0.25 cm}

{\bf Lemma 1.} {\it Let $F$ be a closed set of measure zero on $T$. Then $F$ can be covered by a finite number of open 
intervals of $T$ such that the sum of their lengths is arbitrarily small.}

\vspace{0.25 cm}

The next analogous lemma does not need in a proof either.

 \vspace{0.25 cm} 
 
{\bf Lemma 2.} {\it The set $F$ can be covered by an open set $G_n$ consisting of a finite number of disjoint open 
intervals of $T$ such that $m(G_n) < \frac{1}{2^n}$, $\overline G_{n+1} \subset G_n$, and $F=\bigcap_{n=1}^\infty G_n$.}

\vspace{0.25 cm}

Denote by  $I_m^{(n)}$,   $m = 1, 2, ..., k_n,$ the disjoint open intervals of $T$ constituting the set $G_n$; 
thus\footnote{If $1 \notin F$, then we assume $1 \notin G_n$ for all $n$ (this will simplify the proofs).}
$G_n=\bigcup_{m=1}^{k_n}I_m^{(n)}$. 
Of course, we assume that $I_m^{(n)} \cap F \neq \emptyset$ for each $m = 1, 2, ..., k_n$. 

Let $J_m^{(n)}$ be an open subinterval of $I_m^{(n)}$, such that  $\overline {J_m^{(n)}} \subset I_m^{(n)}$ and $J_m^{(n)} \cap F  = I_m^{(n)} \cap F $  ($m = 1, 2, ..., k_n$). 
Denote $$E_n=\bigcup_{m=1}^{k_n} \overline {J_m^{(n)}}$$ for each natural $n$.
Obviously $F \subset E_n \subset G_n$ and thus $F=\bigcap_{n=1}^\infty E_n$. (With no loss of generality we may assume $E_{n+1} \subset E_n$;
for example, one can just replace $E_{n+1}$ by  $E_{n+1} \cap E_n$). 

\pagebreak

{\bf Lemma 3.} {\it Let $F$, $G_n$,  $I_m^{(n)}$,  $J_m^{(n)} $, and $E_n$ be as above. 
Then there exists a continuously differentiable and periodic 
function $\phi_n$ on $[0, 2 \pi]$, $0 \leq \phi_n(\theta) \leq 1,$
such that $\phi_n(\theta)=0$ if $e^{i\theta} \in T \setminus G_n$ 
and $\phi_n(\theta)=1$ if $e^{i\theta} \in E_n$.}

\vspace{0.25 cm}

{\bf Proof.} For a fixed $n$, we have only finite number of (disjoint) intervals  $I_m^{(n)}$, and $\overline {J_m^{(n)}} \subset I_m^{(n)}$
($m = 1, 2, ..., k_n$). Set $\gamma_{n,m}(\theta) = 0$ if $e^{i\theta} \in T \setminus I_m^{(n)}$ and $\gamma_{n,m}(\theta) = 1$ if 
$e^{i\theta} \in \overline {J_m^{(n)}}$. 
The graph of
$\gamma_{n,m}$
 consists of three horizontal line segments (on the lines $y=0$ and $y=1$).
Using appropriate pieces of shiftings of cosine function
we ``smoothly" join these segments, thus, extending $\gamma_{n,m}$ up
to a 
periodic and continuously differentiable on $[0, 2 \pi]$ function $\phi_{n,m}$ such that
$0 \leq \phi_{n,m}(\theta) \leq 1$. Denote $\phi_n(\theta)=\sum_{m=1}^{k_n}\phi_{n,m}(\theta)$. The function
$\phi_n$ has all the necessary properties and the proof of Lemma 3 is complete.

\vspace{0.25 cm}

Denote 
\begin{equation}
\phi(\theta)=\sum_{n=1}^{\infty}\phi_n(\theta),
 \end{equation}
which is summable on $[0, 2 \pi]$
since  $\int_0^{2\pi} \phi(\theta) d\theta = \sum_{n=1}^{\infty}\int_0^{2\pi}\phi_n(\theta)d\theta < \sum_{n=1}^{\infty}\frac{1}{2^n}=1.$

\vspace{0.25 cm}

{\bf Lemma 4.} {\it At each $\theta_0$ such that $e^{i\theta_0} \in T \setminus F$ the function $\phi$ is finite and continuously differentiable. }

\vspace{0.25 cm}

Lemma 4 is obvious since, as Lemma 3 implies, at certain neighborhood of $\theta_0$ all but finite number of (continuously differentiable) functions $\phi_n$ vanish. 

Let $u_n(z)$ and $u(z)$ be the Poisson integrals of $\phi_n(\theta)$ and $\phi(\theta)$, respectively. They are positive harmonic functions in $\Delta$.
The function $u_n(z)$, as a Poisson integral of a continuous function, is  continuous on $\overline \Delta$, and $u_n(e^{i\theta})=\phi_n(\theta)$ for all $e^{i\theta} \in T $. 
By Lemma 4 the function $\phi$ is  
continuous at each $\theta_0$ such that
 $e^{i\theta_0} \in T \setminus F$ and the elementary 
property of Poisson integrals immediately implies the following lemma.

\vspace{0.25 cm}

{\bf Lemma 5.} {\it Define the function $u(z)$ at each point $e^{i\theta_0} \in T \setminus F$ by setting $u(e^{i\theta_0})=\phi(\theta_0).$ Then the (extended) function $u(z)$ is continuous on the set $\overline \Delta \setminus F$.} 

\vspace{0.25 cm}

Thus  $u(z)$ is defined and continuous on $\overline \Delta \setminus F$. We also define $u(z)=\infty$ on  $F.$ 
 If $l$ is any natural number and $z \in\overline\Delta$, we have   \begin{equation}
u(z) \geq \sum_{n=1}^{l}u_n(z).\end{equation}   
Indeed, $\phi_n(\theta) \geq 0$ implies  (2) for $z \in \Delta$. If $z \in F$, (2) is evident as $u(z) = \infty$ on $F$. Finally, if $z \in T \setminus F,$ then
$u(z)=  u(e^{i\theta}) = \phi(\theta)  \geq \sum_{n=1}^{l}\phi_n(\theta) = \sum_{n=1}^{l}u_n(e^{i\theta}) =\sum_{n=1}^{l}u_n(z)$.

\vspace{0.25 cm}

{\bf Lemma 6.} {\it Let $u(z)$ be defined on the closed unit disk  $\overline \Delta$ as above and let $e^{i\theta_0} \in F$. 
Then $u(z) \to \infty$ as  $z$ ($z \in \overline \Delta$) tends to $e^{i\theta_0}$ arbitrarily.}

\vspace{0.25 cm}

{\bf Proof.} Since $e^{i\theta_0} \in F$,  we have $e^{i\theta_0} \in E_n$ for each $n$, and, by Lemma 3, $u_n(e^{i\theta_0}) = \phi_n(\theta_0)=1.$ Thus 
$\sum_{n=1}^{l}u_n(e^{i\theta_0}) =l.$ Since $l$ is arbitrarily large and each $u_n(z)$ is continuous on $\overline \Delta$, (2) implies that  $u(z) \to \infty$ when
$z \to e^{i\theta_0}$ ($z \in \overline \Delta$). Lemma 6 is proved. 

\vspace{0.25 cm}

{\bf Proof of Theorem A.}
Let $v(z)$ be a conjugate harmonic function of $u(z).$ By Lemma 4 the function $\phi$ is continuously differentiable  on $T \setminus F$.
Thus, by the well known property of the conjugate function (see, e.g. \cite{Hoff}, p. 79), 
$v(z)$ can be extended continuously to the set $T \setminus F.$ Thus, $v(z)$ is continuous on 
$\overline \Delta \setminus F.$ 
Taking into account also Lemma 5 and Lemma 6, we conclude that the function $\omega_F(z) =\frac{1}{1+u(z)+iv(z)}$ belongs to the disk 
algebra, and $\omega_F(z) =0$ if and only if $z \in F.$ Theorem A is proved. 

\vspace{0.25 cm}

To prove Theorem A$'$, note that the function $\lambda_F (z) =\frac{u(z)+iv(z)}{1+u(z)+iv(z)}$ 
belongs to the disk algebra and has the properties formulated in Theorem A$'$.

\vspace{0.25 cm}

{\bf Remark.} The above proof of Lemma 6 shows that $u(z)$ approaches to $\infty$ {\it uniformly} as $z$ ($z \in \overline \Delta$) approaches to any point of the set $F.$

\vspace{0.25 cm}

\end{section}

\begin{section}{Some further remarks}

An important ingredient of Fatou's proof of Theorem A is the following proposition: 

\vspace{0.25 cm}

{\it If the series with positive terms $\sum_{k=1}^{\infty} l_k$ converges, then 
there exists a sequence $\{A_k\},$ with $ A_k > 0, \  A_k \rightarrow +\infty,$
such that $\sum_{k=1}^{\infty} l_kA_k$ converges.} 

\vspace{0.25 cm}

The above proof does not need in this proposition.
Another advantage of our proof is that it readily implies that the harmonic function $u(z)$ tends 
to $\infty$ as $z$ ($z \in \overline \Delta$) approaches to the points of $F.$
In Fatou's proof the same conclusion does not follow so easily and
a reference to 
the theorem on the infinite (unrestricted) boundary values is needed.

\vspace{0.25 cm}

Using the above approach, we derive an elementary (self contained) proof for Lemma 5.

By (1) for $z \in \Delta$ we have 

$$u(z)= \sum_{k=1}^{\infty} u_k(z) = \sum_{k=1}^n u_k(z) + \sum_{k=n+1}^{\infty} u_k(z)   
    =  W_n(z) + \tilde W_n(z),$$
where $W_n(z)$ and $\tilde W_n(z)$ denote the partial sum and the reminder series, respectively. 
Fix an arbitrary point $e^{i\theta_0} \in T \setminus F.$ By Lemma 3, there exists a  $\delta>0$ 
such that all but finite number of functions $\phi_n$ vanish on
the neighborhood 
 $B_\delta(\theta_0)=(\theta_0 -\delta,\theta_0+\delta).$
Fix $N$ so large that if $k>N$
then
$\phi_k$ vanishes on $B_\delta(\theta_0),$ and denote
$\tilde \phi_N (\theta) = \sum_{k=N+1}^{\infty}\phi_k(\theta).$
Then $\tilde W_N(z)$ is the Poisson integral of the
summable function $\tilde \phi_N (\theta)$, which vanishes on $B_\delta(\theta_0).$
Thus $$\tilde W_N(z) \equiv \tilde W_N(re^{i\theta})  =\frac{1}{2\pi} \int_{-\pi}^{\pi}\tilde \phi_N(t)P_r(\theta - t)dt =\frac{1}{2\pi} \int_A \tilde \phi_N(t)P_r(\theta - t) dt,$$
where $A = [-\pi, \pi] \setminus B_\delta(\theta_0).$ 
If $t \in A$ and $\theta \in B_{\frac{1}{2}\delta}(\theta_0) = (\theta_0 - \frac{1}{2} \delta,\theta_0+ \frac{1}{2} \delta)$, then
$|\theta - t| \geq \frac{\delta}{2}.$ Thus, for $\theta \in B_{\frac{1}{2}\delta}(\theta_0),$
$$|\tilde W_N(re^{i\theta})| =  \frac{1}{2\pi} \int_A \tilde \phi_N(t)P_r(\theta - t) dt
 \leq \max_{\frac{\delta}{2} \leq |\theta - t|} \{P_r(\theta -t)\} \frac{1}{2\pi} \int_A \tilde \phi_N(t) = \max_{\frac{\delta}{2} \leq |t|} \{ P_r(t)\} \frac{1}{2\pi} \int_A \tilde \phi_N(t).$$
 The last quantity tends to zero as $r \rightarrow 1.$ Thus, $\tilde W_N(re^{i\theta})$ can be extended
 continuously at the points of the set $B_{\frac{1}{2}\delta}(\theta_0)$, where its values are zero.  In particular, $\tilde W_N(z)$ is continuous at  $e^{i\theta_0}.$
 Since $W_N(z)$ is continuous on $|z| \leq 1$, $u(z)$ is continuous at $e^{i\theta_0}.$ Lemma 5 is proved.

\vspace{0.25 cm}

In closing we note that some recent applications of Fatou's interpolation theorem are given in \cite{Dan1} and \cite{Dan2}.

\end{section}

\begin{minipage}[t]{6.5cm}
Arthur A. Danielyan\\
Department of Mathematics and Statistics\\
University of South Florida\\
Tampa, Florida 33620\\
USA\\
{\small e-mail: adaniely@usf.edu}
\end{minipage}


\begin{thebibliography}{99}


\bibitem{Dan1} A.A. Danielyan, A theorem of Lohwater and Piranian, Proc. of AMS, {\bf 144}, (2016), 3919-3920.

\bibitem{Dan2} A.A. Danielyan, Fatou's interpolation theorem implies the Rudin-Carleson theorem, J. Fourier
Anal. Appl., {\bf 23}, (2017), 656-659.

\bibitem{Fa} P. Fatou, S\'{e}ries trigonom\'{e}triques et s\'{e}ries de Taylor, Acta Math., {\bf 30}, (1906), 335-400.

\bibitem{Gar} J.B. Garnett, Bounded Analytic Functions, Academic
Press, 1981.

\bibitem{Hoff} K. Hoffman, Banach Spaces of Analytic Functions,
Prentice Hall, Englewood Cliffs, New Jersey, 1962.

\bibitem{Koos} P. Koosis, Introduction to $H^p$ Spaces, Cambridge University Press, Cambridge, 1998.


\bibitem{Priv} I.I. Privalov, Boundary Properties of Analytic Functions, Second ed., GITTL,
Moscow-Leningrad, 1950 (in Russian).


\end{thebibliography}
\end{document}